\documentclass[10pt]{amsart} 

\usepackage[utf8]{inputenc} 

\usepackage{geometry} 
\usepackage[dvipdfmx]{graphicx} 

\usepackage[usenames]{xcolor}
\usepackage{amssymb}
\usepackage{amsmath}

\newtheorem{theorem}{Theorem}

\newtheorem{proposition}[theorem]{Proposition}

\newtheorem{definition}[theorem]{Definition\rm}

\newcommand{\op}[1]{{\operatorname{#1}}}
\renewcommand{\L}{\op{L}}
\renewcommand{\exp}{\op{exp}}
\newcommand{\tder}{\op{tder}}
\newcommand{\KV}{\op{KV}}
\newcommand{\Sol}{\op{Sol}}
\renewcommand{\ell}{\text{ell}}
\renewcommand{\div}{\op{div}}

\newcommand{\dop}[1]{%
	\expandafter\newcommand\csname #1\endcsname{{\op{#1}}}
}

\dop{TAut}
\dop{tr}
\dop{ad}
\dop{gr}
\dop{Der}

\newcommand{\gn}{{(g,n+1)}}
\newcommand{\g}{\mathfrak{g}}

\newcommand{\Q}{\mathbb Q}

\begin{document}
\title{Higher genus Kashiwara-Vergne problems and the Goldman-Turaev Lie bialgebra}
\author{Anton Alekseev}
\address{Department of Mathematics, University of Geneva, 2-4 rue du Lievre, 1211 Geneva, Switzerland}
\email{anton.alekseev@unige.ch}

\author{Nariya Kawazumi}
\address{Department of Mathematical Sciences, University of Tokyo, 3-8-1 Komaba, Meguro-ku, Tokyo 153-8914, Japan}
\email{kawazumi@ms.u-tokyo.ac.jp}

\author{Yusuke Kuno}
\address{Department of Mathematics, Tsuda College, 2-1-1 Tsuda-machi, Kodaira-shi, Tokyo 187-8577, Japan}
\email{kunotti@tsuda.ac.jp}

\author{Florian Naef}
\address{Department of Mathematics, University of Geneva, 2-4 rue du Lievre, 1211 Geneva, Switzerland}
\email{florian.naef@unige.ch}

\maketitle

\begin{abstract}
We define a family $\KV^\gn$ of Kashiwara-Vergne problems associated with compact connected oriented 2-manifolds of genus $g$ with $n+1$ boundary components. The problem $\KV^{(0,3)}$ is the classical Kashiwara-Vergne problem from Lie theory. We show the existence of solutions of 
$\KV^\gn$ for arbitrary $g$ and $n$. The key point is the solution of $\KV^{(1,1)}$ based on the results by B. Enriquez on elliptic associators. Our construction is motivated by applications to the formality problem for the Goldman-Turaev Lie bialgebra $\g^{(g, n+1)}$. In more detail, we show that every solution of $\KV^\gn$ induces a Lie bialgebra isomorphism between  $\g^{(g, n+1)}$  and its associated graded $\op{gr} \, \g^{(g, n+1)}$. For $g=0$, a similar result was obtained by G. Massuyeau using the Kontsevich integral.

This paper is a summary of our results. Details and proofs will appear elsewhere.
\end{abstract}

\section{The Goldman-Turaev Lie bialgebra}
Let $\Sigma = \Sigma_{g, n+1}$ be an oriented surface of genus $g$ with $n+1$ boundary components. We fix a framing (that is, a trivialization of the tangent bundle) of $\Sigma$ and choose a base point $* \in \partial \Sigma$. We denote by $\pi = \pi_1(\Sigma, *)$ the fundamental group of $\Sigma$ and  define $\g^{(g,n+1)} = \Q[S^1, \Sigma]  =  \Q \pi / [\Q \pi, \Q \pi]$ to be the vector space spanned by  homotopy classes of free loops. When no confusion arises, we shorten the notation to $\g$.

The vector space $\g^{(g,n+1)}$ carries a canonical Lie bialgebra structure defined in terms of intersections of loops.
The Lie bracket on $\g^{(g,n+1)}$ is called the {\em Goldman bracket} \cite{Go86} and is defined as follows.
Let $\alpha$ and $\beta$ be loops on $\Sigma$ whose intersections are transverse double points.
Then, the Lie bracket $[\alpha,\beta]$ is given by
\[
[\alpha,\beta]=\sum_{p\in \alpha \cap \beta} \varepsilon_p\, \alpha *_p \beta,
\]
where $\varepsilon_p\in \{ \pm 1\}$ is the local intersection number of $\alpha$ and $\beta$ at $p$, and $\alpha *_p \beta$ is the homotopy class of the concatenation of the loops $\alpha$ and $\beta$ based at $p$ (see Fig. 1).
\begin{figure}
\centering
\input{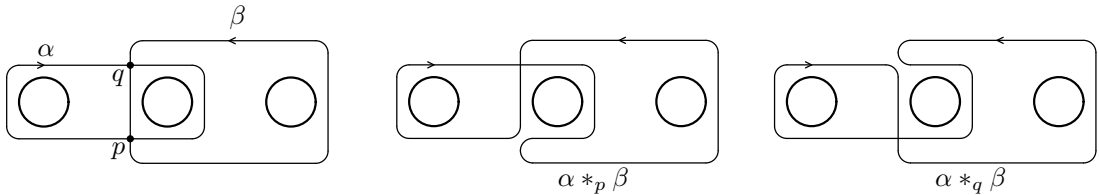}
\centering
\caption{Goldman bracket.
In this figure, $[\alpha,\beta]=\alpha *_p \beta-\alpha *_q \beta$.}
\label{fig:bracket}
\end{figure}
The Lie cobracket on $\g^{(g,n+1)}$ is called the {\em Turaev cobracket} \cite{Tu91} and is defined as follows.
Let $\gamma$ be a loop on $\Sigma$.
By a suitable homotopy, we can deform $\gamma$ to an immersion with transverse double points whose rotation number with respect to the framing of $\Sigma$ is zero.
For each self-intersection $p$ of $\gamma$, one can divide $\gamma$ into two branches $\gamma^1_p$ and $\gamma^2_p$,
where the pair of the tangent vectors of $\gamma^1_p$ and $\gamma^2_p$ forms a positive basis for $T_p\Sigma$.
Then, the Lie cobracket $\delta(\gamma)$ is given by 
\[
\delta(\gamma)=\sum_p \gamma^1_p\otimes \gamma^2_p-\gamma^2_p\otimes \gamma^1_p,
\]
where the sum is taken over all the self-intersections of $\gamma$ (see Fig. 2).
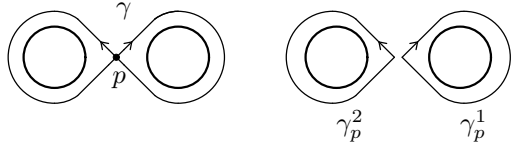
\begin{figure}
\centering
{\unitlength 0.1in%
\begin{picture}(26.0000,6.4000)(1.2000,-7.5000)%
%
\special{pn 13}%
\special{ar 360 440 160 160 0.0000000 6.2831853}%
%
\special{pn 13}%
\special{ar 1000 440 160 160 0.0000000 6.2831853}%
%
\special{pn 8}%
\special{pa 640 480}%
\special{pa 720 400}%
\special{fp}%
%
\special{pn 8}%
\special{pa 640 400}%
\special{pa 720 480}%
\special{fp}%
%
\special{pn 8}%
\special{pa 640 480}%
\special{pa 480 640}%
\special{fp}%
%
\special{pn 8}%
\special{ar 360 480 200 200 0.9272952 1.5707963}%
%
\special{pn 8}%
\special{ar 360 440 240 240 1.5707963 4.7123890}%
%
\special{pn 8}%
\special{ar 360 400 200 200 4.7123890 5.3558901}%
%
\special{pn 8}%
\special{pa 640 400}%
\special{pa 480 240}%
\special{fp}%
%
\special{pn 8}%
\special{ar 1000 440 240 240 4.7123890 1.5707963}%
%
\special{pn 8}%
\special{ar 1000 400 200 200 4.0688879 4.7123890}%
%
\special{pn 8}%
\special{pa 720 400}%
\special{pa 880 240}%
\special{fp}%
%
\special{pn 8}%
\special{ar 1000 400 200 200 4.0688879 4.7123890}%
%
\special{pn 8}%
\special{ar 1000 480 200 200 1.5707963 2.2142974}%
%
\special{pn 8}%
\special{pa 880 640}%
\special{pa 720 480}%
\special{fp}%
%
\special{pn 8}%
\special{pa 760 360}%
\special{pa 752 400}%
\special{fp}%
\special{pa 760 360}%
\special{pa 720 368}%
\special{fp}%
%
\special{pn 8}%
\special{pa 600 360}%
\special{pa 608 400}%
\special{fp}%
\special{pa 600 360}%
\special{pa 640 368}%
\special{fp}%
%
\special{pn 13}%
\special{ar 1820 440 160 160 0.0000000 6.2831853}%
%
\special{pn 13}%
\special{ar 2460 440 160 160 0.0000000 6.2831853}%
%
\special{pn 8}%
\special{ar 1800 480 200 200 0.9272952 1.5707963}%
%
\special{pn 8}%
\special{ar 1800 440 240 240 1.5707963 4.7123890}%
%
\special{pn 8}%
\special{ar 1800 400 200 200 4.7123890 5.3558901}%
%
\special{pn 8}%
\special{pa 2040 360}%
\special{pa 2048 400}%
\special{fp}%
\special{pa 2040 360}%
\special{pa 2080 368}%
\special{fp}%
%
\special{pn 8}%
\special{pa 1920 240}%
\special{pa 2120 440}%
\special{fp}%
%
\special{pn 8}%
\special{pa 1920 640}%
\special{pa 2120 440}%
\special{fp}%
%
\special{pn 8}%
\special{ar 2480 480 200 200 1.5707963 2.2142974}%
%
\special{pn 8}%
\special{ar 2480 440 240 240 4.7123890 1.5707963}%
%
\special{pn 8}%
\special{ar 2480 400 200 200 4.0688879 4.7123890}%
%
\special{pn 8}%
\special{pa 2240 360}%
\special{pa 2232 400}%
\special{fp}%
\special{pa 2240 360}%
\special{pa 2200 368}%
\special{fp}%
%
\special{pn 8}%
\special{pa 2360 240}%
\special{pa 2160 440}%
\special{fp}%
%
\special{pn 8}%
\special{pa 2360 640}%
\special{pa 2160 440}%
\special{fp}%
\put(6.6000,-6.0000){\makebox(0,0)[lb]{$p$}}%
\put(6.8000,-2.4000){\makebox(0,0)[lb]{$\gamma$}}%
\put(24.6000,-8.8000){\makebox(0,0)[lb]{$\gamma_p^1$}}%
\put(18.2000,-8.8000){\makebox(0,0)[lb]{$\gamma_p^2$}}%
%
\special{pn 4}%
\special{sh 1}%
\special{ar 680 440 16 16 0 6.2831853}%
\end{picture}}%
\caption{Turaev cobracket.
In this figure, $\delta(\gamma)=\gamma_p^1\otimes \gamma_p^2-\gamma_p^2\otimes \gamma_p^1$.}
\label{fig:cobracket}
\end{figure}

The group algebra $\Q \pi$ carries a canonical filtration with the following property. Choose a set of generators 
$\alpha_i, \beta_i, \gamma_j \in \pi$ with $i=1, \dots, g, j = 1, \dots, n$ such that 
$$
\prod_{i=1}^g [\alpha_i, \beta_i] \prod_{j=1}^n \gamma_j = \gamma_0,
$$
where $\gamma_0$ is the homotopy class of the boundary component which contains $*$.
Then, the elements $(\alpha_i-1), (\beta_i -1) \in \Q \pi$ are of filtration degree 1 and $(\gamma_j-1) \in \Q\pi$ have filtration degree 2.  This filtration on $\pi$ induces a two step filtration on $H=H_1(\Sigma, \Q)$ with $H^{(1)} = H$ and $H^{(2)}$ the kernel of the intersection pairing. 
The associated graded of $\widehat{\Q \pi}$ with respect to this filtration is the complete Hopf algebra  given by $\op{gr}\, \widehat{\Q \pi} \cong T(\op{gr}\, H) = U(\L(\op{gr}\, H))$, where $T(\op{gr}\, H)$ is the completed tensor algebra of the graded vector space $\op{gr}\, H$ and $U(\L(\op{gr}\, H))$ is the completed universal enveloping algebra of the free Lie algebra $\L(\op{gr}\, H)$.

Let $\hat{\g}^{(g,n+1)}$ be the completion of $\g^{(g, n+1)}$. The completed associated graded vector space 
$\op{gr}\, \g$ can be canonically identified with the space of formal series in cyclic tensor powers of the graded vector space $\op{gr} \, H = H/H^{(2)} \oplus H^{(2)}$. A choice of a basis in $\pi$ induces a basis in $\op{gr} \, H: x_i, y_i, z_j $ with $i=1, \dots, g, j=1, \dots n$, where 
$x_i, y_i$ are of degree 1 and $z_j$ are of degree 2. 
The graded vector space $\op{gr} \, \g$ carries a canonical Lie bialgebra structure induced by the Goldman-Turaev Lie bialgebra structure on $\g$. It turns out that both the Lie bracket and Lie cobracket on $\op{gr}\, \g$ are of degree $(-2)$.

Let $\mathcal{M}(\Sigma)$ be the mapping class group of the surface $\Sigma$ fixing the boundary $\partial \Sigma$ pointwise.
There are a subset $\mathcal{M}(\Sigma)^{\circ}\subset \mathcal{M}(\Sigma)$ and an embedding $\tau: \mathcal{M}(\Sigma)^{\circ} \to \hat{\g}^{\gn}$ such that any Dehn twist is in $\mathcal{M}(\Sigma)^{\circ}$, $\mathcal{M}(\Sigma_{g,1})$ includes the Torelli group and the graded quotient of $\tau$ is the classical Johnson homomorphism \cite{KK16}.
In \cite{KK15}, the second and third authors proved that $\delta\circ \tau=0: \mathcal{M}(\Sigma)^{\circ}\to \hat{\g}\hat{\otimes}\hat{\g}$.
This is one of the motivations to study the graded version of the Goldman-Turaev Lie bialgebra and the corresponding formality problem.

\begin{definition}
A {\em group-like expansion} is an isomorphism $\theta: \widehat{\Q \pi} \to  \op{gr} \, \widehat{\Q \pi}$ of complete filtered Hopf algebras with the property $\op{gr} \, \theta = {\rm Id}$.
\end{definition}

It is easy to see that group-like expansions exist. For instance, every choice of a basis in $\pi$ described above induces a group-like expansion $\theta^{\op{exp}}$ defined by its values on generators:
$$
\theta^{\op{exp}}(\alpha_i) = e^{x_i}, \quad
\theta^{\op{exp}}(\beta_i) = e^{y_i}, \quad
\theta^{\op{exp}}(\gamma_j) = e^{z_j}.
$$
In fact, group-like expansions are a torsor under the group of automorphisms of the complete Hopf algebra $\op{gr}\, \widehat{\Q \pi} \cong T(\op{gr} \, H)$ with associated graded the identity. That is, every group-like expansion is of the form
$\theta=F \circ \theta^\exp$ for some $F \in \op{Aut} (\L(\op{gr} \, H))$. Furthermore, every group-like expansion defines an isomorphism of filtered vector spaces $\g \to \op{gr} \, \g$ with associated graded the identity map.

\begin{definition}
A group-like expansion $\theta$ is called {\em homomorphic} if it induces an isomorphism of Lie bialgebras 
$\g \to \op{gr} \, \g$.
\end{definition}

It is easy to check that expansions $\theta^{\op{exp}}$ are not homomorphic.
Existence of homomorphic expansions is one of the main results of this paper. Our strategy is to reformulate the problem in terms of properties of the automorphism $F \in \op{Aut} (\L(\op{gr} \, H))$. This leads us to a generalization of the Kashiwara-Vergne problem in the theory of free Lie algebras.

\section{The Kashiwara-Vergne problem in higher genus}
Denote by $\L^{(g,n+1)} = \L(\op{gr} \, H) $ the completed free Lie algebra in the generators $x_i, y_i, z_j$ with $\deg x_i = \deg y_i =1, \deg z_j =2$. Define the completed graded Lie algebra of {\em tangential derivations} (in the sense of \cite{AT}):
$$
\tder^{(g, n+1)} = \{ u \in \Der^+ (\L^{(g,n+1)}) \ | \ u(z_j) = [z_j, u_j] \  \text{for some $u_j \in \L^{(g,n+1)} $} \} ,
$$
This Lie algebra integrates to a pro-unipotent group
$$
\op{TAut}^{(g,n+1)} = \{ F \in \op{Aut}^+(\L^{(g,n+1)}) \ | \ F(z_j) = F^{-1}_j z_j F_j \ \text{for some $F_j \in \exp(\L^{(g,n+1)})$} \} .
$$
 Also following \cite{AT}, let $\op{tr}^{(g,n+1)}$ denote the vector space of series in cyclic words in the $x_i, y_i$ and $z_j$'s and $\op{tr}$ the natural projection from associative words to cyclic words.
The space $\op{tr}^{(g, n+1)}$ carries an action of $\tder^{\gn}$ coming from a natural action of $\Der^+ (\L^{\gn})$.
Recall the definition of the map $\partial_{x_i} : \L^{(g,n+1)} \to U(\L^{(g,n+1)})$ given by the formula
\begin{align*}
\alpha(x_1, \cdots, x_i + \epsilon \xi, \cdots, x_g, y_1, \cdots, y_g, z_1, \cdots z_n) &= \alpha + \epsilon \, \op{ad}(\partial_{x_i} \alpha) \,  \xi + O(\epsilon^2) \quad \text{for $\alpha \in \L^{(g, n+1)}$},
\end{align*}
and the same definition for $y_i$ and $z_i$.
One defines the non-commutative  divergence map $\op{div} : \op{tder}^{(g,n+1)} \rightarrow \op{tr}^{(g,n+1)} $ as follows:
$$
u \longmapsto \sum_{i=1}^{g} \op{tr}(\partial_{x_i}(u(x_i)) + \partial_{y_i}(u(y_i))) + \sum_{j=1}^n \op{tr}(z_j \partial_{z_j}(u_j)),
$$
where $u_j$ is chosen such that it has no linear terms in $z_j$.
In the case of $n=0$, the divergence coincides with  the Enomoto-Satoh obstruction \cite{ES} for  surjectivity of the Johnson homomorphism.

\begin{proposition}
The map $\op{div}$ is a Lie algebra 1-cocycle on $ \op{tder}^{(g,n+1)}$ with values in $\op{tr}^{(g,n+1)}$. It integrates to a group 1-cocycle $j: \op{TAut}^{(g,n+1)}  \to \op{tr}^{(g,n+1)}$.
\end{proposition}


Let $r(s)$ be the power series $\op{log}(\frac{s}{e^s-1}) \in \Q[[s]]$, and let $\mathbf{r} = \sum_{i=1}^g  \op{tr}(r(x_i) + r(y_i) ) \in \op{tr}^{(g,n+1)}$. Below we define a new family of {\em Kashiwara-Vergne problems} associated to a surface of genus $g$ with $n+1$ boundary components.

\begin{definition}[KV Problem of type $(g, n+1)$]
Find an element $F \in \op{TAut}^{(g,n+1)}$ such that
\begin{gather*}
F(\sum_{i=1}^{g} [x_i, y_i] + \sum_{j=1}^n z_j) = 
\op{log} ( \prod_{i=1}^{g} (e^{x_i} e^{y_i} e^{-x_i} e^{-y_i})  \, \prod_{j=1}^{n} e^{z_j} ) =: \xi  
\tag{KVI$^\gn$} \\
j(F) = \sum_{i=1}^n \op{tr} \,  h(z_i)  - \op{tr} \, h ( \xi)  - \mathbf{r} \quad \text{for some Duflo function $h \in \Q[[s]]$} \tag{KVII$^\gn$}.
\end{gather*}

\end{definition}

Let $\Sol \KV^{(g,n+1)}$ denote the set of solutions of the KV problem of type $(g, n+1)$.
Note that $\KV^{(0,3)}$ is the classical KV problem  \cite{KV} in the formulation of \cite{AT}.
The following is the first main result of this note:

\begin{theorem}  \label{thm:main1}
Let $F\in \Sol \KV^{(g,n+1)}$.
Then, the group-like expansion $F^{-1} \circ \theta^\exp$ is homomorphic.
Moreover, if $F\in \op{TAut}^{(g,n+1)}$ satisfies $(\KV {\rm I}^{\gn})$, then
$F^{-1}\circ \theta^\exp$ is homomorphic if and only if 
$F\in \Sol \KV^{(g,n+1)}$.
\end{theorem}

This result holds true for the framing {\em adapted} to the set of generators $\alpha_i, \beta_i, \gamma_j \in \pi$.
Note that any framing on $\Sigma$ is specified by the values of its rotation number function on simple closed curves which are freely homotopic to $\alpha_i,\beta_i,\gamma_j$.
The adapted framing takes values $0$ on $\alpha_i,\beta_i$, and $-1$ on $\gamma_j$.
Other framings require a more detailed discussion.
Our proof of Theorem \ref{thm:main1} uses the theory of van den Bergh double brackets \cite{vdB08}, their moment maps \cite{N15} and their relation to the Goldman bracket \cite{MT13}.

\section{Solving $\KV^{\gn}$}
Our second main result is  the following theorem:

\begin{theorem} \label{thm:main2}
For all $g \geq 0, n\geq0$, $ \Sol \KV^{(g,n+1)}\neq \varnothing$.
\end{theorem}

Together, Theorem \ref{thm:main1} and Theorem \ref{thm:main2} imply existence of homomorphic expansions for any $g$ and $n$ (for $g=0$, an independent proof was given by G. Massuyeau \cite{Mas15}). Among other things, it follows that the obstruction to surjectivity of the Johnson homomorphism provided by the Turaev cobracket is equivalent to the Enomoto-Satoh obstruction. 

In what follows, we sketch a proof of this statement. 
%
%
Let us denote the variables appearing in the definition of $\tder^{(g_1+g_2, n_1 + n_2+1)}$ by $x^1_i,y^1_i,z^1_j$ and $x^2_i,y^2_i,z^2_j$ respectively and define the following map
\begin{align*}
\mathcal{P}: \tder^{(0,3)} &\longrightarrow \tder^{(g_1+g_2, n_1+n_2+1)} \\
(u_1,u_2) & \longmapsto \left( w^k \mapsto [w^k, u_k(\phi_1,\phi_2)] \right), \text{with $w^k\in \{x_i^k,y_i^k,z_j^k\}, k=1,2$},
\end{align*}
where
$$ 
\phi_1 = \sum_i [x^1_i,y^1_i] + \sum_j z^1_j,\quad
\phi_2 = \sum_i [x^2_i,y^2_i] + \sum_j z^2_j.
$$
This map is a Lie algebra homomorphism. It lifts to a group homomorphism $\op{TAut}^{(0,3)} \to \op{TAut}^{(g_1+g_2,n_1+n_2+1)}$ (also denoted by $\mathcal{P}$).
Denote by $t \in \tder^{(0,3)}$ the tangential derivation $t: z_1 \mapsto [z_1, z_2], z_2 \mapsto [z_2, z_1]$. Recall that for $F \in \Sol \KV^{0,3}$ there is a family of solutions $F_\lambda=F \exp(\lambda t)$ for $\lambda \in \Q$.

\begin{proposition} \label{prop:glueing}
Let $F_1 \in \Sol \KV^{(g_1, n_1+1)}, F_2 \in \Sol \KV^{(g_2, n_2+1)}, F \in \Sol \KV^{(0, 3)}$ such that their Duflo functions coincide, $h_1=h_2=h$.
If $n_1=n_2=0$ or $g_1=0$ or $g_2=0$, then there is $\lambda \in \Q$ such that
$$
\tilde{F} := (F_1 \times F_2) \circ \mathcal{P}(F_\lambda) \in \Sol \KV^{(g_1 + g_2, n_1+n_2 +1)},
$$
and the corresponding Duflo function $\tilde{h} = h_1 = h_2 = h$.
\end{proposition}

Proposition \ref{prop:glueing} reduces the proof of Theorem \ref{thm:main2} to the cases $\KV^{(0,3)}$ and $\KV^{(1,1)}$. By \cite{AT}, the problem $\KV^{(0,3)}$ does admit solutions.
%
The remaining case is thus $(g, n+1) = (1,1)$. Let $\TAut_{z_1-z_2}^{(0,3)} \subset \TAut^{(0,3)}$ denote the subgroup which preserves $z_1 - z_2$ up to quadratic terms. One defines the following group homomorphism 
$\TAut_{z_1-z_2}^{(0,3)} \rightarrow \TAut^{(1,1)}$:
$$
F \longmapsto F^\ell: \left\{
	\begin{array}{lr}
		e^{x_1} \mapsto F_1(\psi_1, \psi_2)^{-1} e^{x_1} F_2(\psi_1,\psi_2) \\
		e^{y_1} \mapsto F_2(\psi_1, \psi_2)^{-1} e^{y_1} F_2(\psi_1,\psi_2),
	\end{array}
\right.
\quad
\text{where $\psi_1=e^{x_1}y_1e^{-x_1}$ and $\psi_2=-y_1$.}
$$
Furthermore, let $\varphi \in \TAut^{(1,1)}$ be an automorphism defined by $\varphi(x_1) = x_1,\ \varphi(y_1) = \frac{e^{\ad_{x_1}} -1}{\ad_{x_1}}(y_1)$.

\begin{proposition}   \label{thm:enriquez}
Let $F \in \Sol \KV^{(0,3)}$, then there is a unique $\lambda \in \Q$, such that $(F e^{\lambda t})^{\rm ell} \circ \varphi \in \Sol \KV^{(1,1)}$.
\end{proposition}

The proof of Proposition \ref{thm:enriquez} is based on the results of \cite{Enr}. Together with Proposition \ref{prop:glueing}, it settles in the positive the existence issue for solutions of Kashiwara-Vergne problems $\KV^{(g, n+1)}$. We now turn to the uniqueness problem. Recall the notation
$\phi = \sum_i [x_i, y_i] + \sum_j z_j$.

\begin{definition}
The Kashiwara-Vergne Lie algebra $\mathfrak{krv}^{(g, n+1)}$ is defined as
$$
\mathfrak{krv}^{(g, n+1)} = \{ u \in \tder^{(g, n+1)} \ | \ u(\phi)=0, \div(u) = \sum_j {\rm tr} \, h(z_j)  - {\rm tr} \, h(\phi) \,\, {\rm for} \,\, {\rm some} \,\,
h \in \mathbb{Q}[[s]]\}.
$$
\end{definition}

Of particular interest is the Lie algebra $\mathfrak{krv}^{(1,1)}$:
$$
\mathfrak{krv}^{(1,1)}=\{ u \in \tder^{(1,1)} =\Der^+(\L(x,y)) \ | \  u([x,y])=0, \div(u) = - {\rm tr} \, h([x,y]) \}.
$$
Note that there are other definitions of Kashiwara-Vergne Lie algebras in the literature, see \cite{Sch,SR} for an alternative definition of $\mathfrak{krv}^{(1,1)}$ based on the theory of moulds and \cite{Willwacher} for a graph theoretic definition for arbitrary manifolds. At this point, we do not know what is the relation of these approaches to our considerations.

The pro-nilpotent Lie algebra $\mathfrak{krv}^{(g, n+1)}$ integrates to a group $\mathrm{KRV}^{(g, n+1)}$ which acts freely and transitively on the set of solutions of the Kashiwara-Vergne problem $\KV^{(g, n+1)}$, $G: F \mapsto FG$. 
The following result gives  partial information on the structure of the Lie algebra $\mathfrak{krv}^{(1,1)}$:

\begin{proposition}
\label{prop:krv11}
There is an injective Lie homomorphism of the Grothendieck-Teichm\"uller Lie algebra $\mathfrak{grt}_1$ into
$\mathfrak{krv}^{(1,1)}$
and the elements $\delta_{2n} \in \Der^+ (\L(x,y)), n=1, \dots$ uniquely defined by conditions $\delta_{2n}([x,y])=0, \delta_{2n}(x) = {\rm ad}_x^{2n}(y)$ belong to $\mathfrak{krv}^{(1,1)}$.
\end{proposition}

In \cite{Enr}, it is conjectured that $\mathfrak{grt}_1$ together with $\delta_{2n}$'s form a generating set for the elliptic Grothendieck-Teichm\"uller Lie algebra $\mathfrak{grt}_{\rm ell}$. In view of Proposition \ref{prop:krv11}, we conjecture that $\mathfrak{grt}_{\rm ell}$ injects in $\mathfrak{krv}^{(1,1)}$.

\vskip 0.2cm

{\bf Acknowledgements.} We are indebted to T. Kohno  for a suggestion to start our collaboration. We are grateful to G. Massuyeau, E. Raphael and L. Schneps for sharing with us the results of their work. We thank B. Enriquez, P. Severa and  A. Tsuchiya for fruitful discussions.
The work of A.A. and F.N. was supported in part by the grant 165666 of the Swiss National Science Foundation, by the ERC grant MODFLAT and by the NCCR SwissMAP. Research of N.K. was supported in part by the grants JSPS KAKENHI 15H03617 and 24224002. Y.K. was supported in part by the grant JSPS KAKENHI 26800044.

\end{document}